\documentclass[11pt]{article}
\usepackage{amsfonts, amsmath, amssymb}
\def\qed{\nopagebreak\rule{5pt}{8pt}}
\title{{\bf New Versions of the All-Ones Problem}
\footnote{Supported by NSFC.}}
\author{{ \small Xueliang Li and Xiaoyan Zhang }\\
[3mm]
\small Center for Combinatorics and LPMC\\
\small Nankai University\\
\small Tianjin 300071, P.R. China }

\date{}

\begin{document}

\maketitle{}

\begin{abstract}

We study three new versions of the All-Ones Problem and the
Minimum All-Ones Problem. The original All-Ones Problem is simply
called the Vertex-Vertex Problem, and the three new versions are
called the Vertex-Edge Problem, the Edge-Vertex Problem and the
Edge-Edge Problem, respectively. The Vertex-Vertex Problem has
been studied extensively. For example, existence of solutions and
efficient algorithms for finding solutions were obtained, and the
Minimum Vertex-Vertex Problem for general graphs was shown to be
NP-complete and for trees it can be solved in linear time, etc. In
this paper, for the Vertex-Edge Problem, we show that a graph has
a solution if and only if it is bipartite, and therefore it has
only two possible solutions and optimal solutions. A linear
program version is also given. For the Edge-Vertex Problem, we
show that a graph has a solution if and only if it contains even
number of vertices. By showing that the Minimum Edge-Vertex
Problem can be polynomially transformed into the Minimum Weight
Perfect Matching Problem, we obtain that the Minimum Edge-Vertex
Problem can be solved in polynomial time in general. The Edge-Edge
Problem is reduced to the Vertex-Vertex Problem for the line graph
of a graph. \\
[3mm] {\bf Keywords:} {All-Ones Problem; minimum weight perfect matching;
graph algorithm} \\
[3mm] {\bf MR Subject Classifications(2000):} 05C85, 05C70, 90C27,
68Q25, 68R10\\
\end{abstract}

\section{Introduction}

The term {\it All-Ones Problem} was introduced by Sutner, see
[12]. It has applications in linear cellular automata, see [13]
and the references therein. The problem is cited as follows:
suppose each of the square of an $n\times n$ chessboard is
equipped with an indicator light and a button. If the button of a
square is pressed, the light of that square will change from off
to on and vice versa; the same happens to the lights of all the
edge-adjacent squares. Initially all lights are off. Now, consider
the following questions: is it possible to press a sequence of
buttons in such a way that in the end all lights are on ? This is
referred as the {\it All-Ones Problem}. If there is such a
solution, how to find a such way ? And finally, how to find such a
way that presses as few buttons as possible ? This is referred as
the {\it Minimum All-Ones Problem}. All the above questions can be
asked for arbitrary graphs. Here and in what follows, we consider
connected simple undirected graphs only. One can deal with
disconnected graphs component by component. For all terminology
and notations on graphs, we refer to [2]. An equivalent version of
the All-Ones Problem was proposed by Peled in [10], where it was
called the {\it Lamp Lighting Problem}. The rule of the All-Ones
Problem is called $\sigma^{+}$ rule on graphs, which means that a
button lights not only its neighbors but also its own light. If a
button lights only its neighbors but not its own light, this rule
on graphs is called $\sigma$ rule.

There have been many publications on the All-Ones Problem, see
Sutner [14,15], Barua et al [1] and Dodis and Winkler [4]. Using
linear algebra, Sutner [13] proved that it is always possible to
light every lamp in any graphs by $\sigma^{+}$ rule. Lossers [7]
gave another beautiful proof also by using linear algebra. A
graph-theoretic proof was given by Erikisson et al [6]. In [3] we
gave a graph-theoretic algorithm of cubic time for finding the
solutions. In [11], Sutner proved that the Minimum All-Ones
Problem is NP-complete in general. We gave [3] a linear time
algorithm for finding optimal solutions for trees.\\

In graph-theoretic terminology, a solution to the All-Ones Problem
with $\sigma^{+}$-rule can be stated as follows: given a graph
$G=(V,E)$, where $V$ and $E$ denotes the vertex-set and the
edge-set of $G$, respectively. A subset $X$ of $V$ is a solution
if and only if for every vertex $v$ of $G$ the number of vertices
in $X$ adjacent to or equal to $v$ is odd. Such a subset $X$ is
called an {\it odd parity cover} in [12]. So, the All-Ones Problem
can be formulated as follows: given a graph $G=(V,E)$, does a
subset $X$ of $V$ exist such that for all vertex $v\in V-X$, the
number of vertices in $X$ adjacent to $v$ is odd, while for all
vertex $v\in X$, the number of vertices in $X$ adjacent to $v$ is
even ? If there exists a solution, how to find a one with minimum
cardinality ?  We simply call them the {\it Vertex-Vertex
Problem}, in contrasting with the following new versions of the
All-Ones Problem.\\

From the Vertex-Vertex Problem, one is easily led to propose the
following three problems. These problems have their own interests
and are worth to be studied from algorithmic point of view, even
if not related to the All-Ones Problem. \\

The {\it Vertex-Edge Problem}: suppose there is a light bulb on
each edge and a button on each vertex. If we press the button on a
vertex, the on/off status of the lights on the edges incident with
the vertex will be changed. If the initial status of all lights
are off, we ask whether one can press some buttons to light all
the lights on the edges. This problem can be formulated in
graph-theoretic terminology as follows: given a graph $G=(V,E)$,
does there exist a subset $X$ of $V$ such that for all edge $e\in
E$, exact one of its end-vertices is in $X$ ? If there
exists a solution, how find a one with minimum cardinality ?\\

The {\it Edge-Vertex Problem}: suppose there is a light bulb on
each vertex and a button on each edge. If we press the button on
an edge, the on/off status of the lights on the vertices incident
with the edge will be changed. If the initial status of all lights
are off, we ask whether one can press some buttons to light all
the lights on the vertices. This problem can be formulated in
graph-theoretic terminology as follows: given a graph $G=(V,E)$,
does there exist a subset $X$ of $E$ such that for all vertex
$v\in V$, the number of edges in $X$ incident with $v$ is odd ? If
there exists a solution, how to find a one with minimum
cardinality ? \\

The {\it Edge-Edge Problem}: suppose there is a light bulb and a
button on each edge. If we press the button on an edge, the on/off
status of the light on this edge and those on the edges adjacent
to the edge will be changed. If the initial status of all lights
are off, we ask whether one can press some buttons to light all
the lights on the edges. This problem can also be formulated in
graph-theoretic terminology as follows: given a graph $G=(V,E)$,
does there exist a subset $F$ of $E$ such that for all edge $e\in
E-F$, the number of edges in $F$ adjacent to $e$ is odd, while for
all edge $e\in F$, the number of edges in $F$ adjacent to $e$ is
even ? If there exist a solution, how to find a one with minimum
cardinality ? It is easily seen that this problem is equivalent to
the Vertex-Vertex Problem for the line graph of a
graph.\\

This paper is organized as follows. Section 1 is an introduction.
In Section 2, we show that a graph has a solution to the
Vertex-Edge Problem if and only if it is bipartite, and therefore
it has only two possible solutions and optimal solutions. A linear
programming version is also given. In Section 3, we show that a
graph has a solution to the Edge-Vertex Problem if and only if it
contains even number of vertices. By showing that the Minimum
Edge-Vertex Problem can be polynomially transformed into the
Minimum Weight Perfect Matching Problem, we obtain that the
Minimum Edge-Vertex Problem can be solved in polynomial time in
general. As we already pointed out, in Section 4, the Edge-Edge
Problem is simply reduced to the Vertex-Vertex Problem for the
line graph of a graph, and therefore there always exist solutions
for general graphs, and a solution can be found in polynomial
time.\\

\section{The Vertex-Edge Problem}

The Vertex-Edge Problem and its corresponding minimum problem can
be completely solved, which is the easiest problem among the four
problems. In this section we will find all the graphs that have a
solution for the Vertex-Edge Problem, and figure out the two possible
solutions.\\

\noindent {\bf Theorem 2.1.} The Vertex-Edge Problem has a
solution for a graph $G$ if and only if G is bipartite. Moreover,
a connected bipartite graph $G=(U,V,E)$ has only two solutions $U$
and $V$. Therefore, the solution to the Minimum Vertex-Edge
Problem is $U$ or $V$ such that it attains the value $\min\{|U|,|V|\}$.\\

\noindent {\bf Proof. } For the first statement of the theorem, if
$G$ is bipartite, say $G=(U,V,E)$, then it is obvious that $U$ and
$V$ are solutions to the Vertex-Edge Problem for $G$. Conversely,
if $G$ has a solution for the Vertex-Edge Problem, we only need to
show that there does not exist odd cycle in $G$. Otherwise,
suppose $C$ is an odd cycle in $G$. Then, the lights on the edges
of $C$ should be on in the end. Because all the vertices incident
to the edges of $C$ are on the cycle, there must be a way to press
some of the buttons on the vertices of the cycle $C$ to light all
the lights on the edges of $C$, which is impossible since $C$ is
an odd cycle, a contradiction. Therefore, $G$ does not have odd
cycles, and so $G$ is bipartite.\\

For the second statement of the theorem, suppose that $G=(U,V,E)$
is a connected bipartite graph, and $U_{1}\bigcup V_{2}$ is a
solution to the Vertex-Edge Problem for $G$ such that $U\supseteqq
U_{1}\neq\emptyset, V\supseteqq V_{2}\neq\emptyset$. Consider the
subgraph $(U_{1},V_{1})$ of $G$ induced by $U_{1}$ and the
incident edges, and the subgraph $(U_{2},V_{2})$ of $G$ induced by
$V_{2}$ and the incident edges. We claim that in $G$ there is no
edge between the subgraph $(U_{1},V_{1})$ and the subgraph
$(U_{2},V_{2})$. Otherwise, suppose that there is an edge $e$
between $U_{1}$ and $V_{2}$, then $e$ can not be lighted because
the two end-vertices of $e$ are in the solution. Suppose $e$ is
between $U_{2}$ and $V_{1}$, then $e$ can not be lighted because
both of the two end-vertices of $e$ are not in the solution. In
any case, the edge $e$ can not be lighted by the solution
$U_1\bigcup V_2$, a contradiction. So, $(U_{1},V_{1})$ is not
connected with $(U_{2},V_{2})$, a contradiction to the assumption
that $G$ is connected. Therefore, one of $U_{1}$ and $V_{2}$ must
be empty. If $U_1=\emptyset$, then $V_2=V$ since $U_1\bigcup V_2$
is a solution, and if $V_2=\emptyset$, then $U_1=U$ by the same
argument. Hence, $G$ has only two solutions for the Vertex-Edge Problem. \\

The third statement of the theorem for the Minimum Vertex-Edge
Problem is obvious. The proof is now complete. \qed \\

For a disconnected bipartite graph, we can find the solutions for
its every connected component, and then the union of them is the
solutions for the whole graph.\\

Let $G=(U,V,E)$ be a bipartite graph such that $|U \bigcup V|=n$.
For every vertex of $G$, we assign a variable $x_{i}$. Then, the
Minimum Vertex-Edge Problem can be formulated by an $(0,1)$-linear
program as follows:
\[ \min \sum_{i=1}^{n} x_{i} \]
\[  \left\lbrace
  \begin{array}{l l}
    A^{T}X=I & \\
    x_i\in \{0,1\}\\
\end{array}
\right. \] \\
where $A$ is the edge(row)-vertex(column) incidence matrix of $G$,
$X^{T}=(x_{1},x_{2},\cdots,x_{n})$ and $I^{T}=(1,1,\cdots,1)$.
Because $G$ is an undirected bipartite graph, $A$ must be a
totally unimodular matrix. So, the above programming can be solved
simply by solving its relaxed linear programming, which can be
done in polynomial time. So, the Vertex-Edge Problem and its
corresponding minimum problem can be solved easily.\\

\section{The Edge-Vertex Problem}

Unlike the Vertex-Edge Problem, the Edge-Vertex Problem is not
that easy to solve. At first, it is easy to see that the
Edge-Vertex Problem has the following relation with the ODD SET
Problem [5]: for a graph $G=(V,E)$, if we take the set of its
edges as the set $R=E=\{e_{1},e_{2},\cdots,e_{m} \}$, the set of
its vertices as the set $B=V=\{v_{1},v_{2},\cdots,v_{n} \}$ and
$e_{i}$ is adjacent to $v_{j}$ if and only if $v_{j}$ is an
end-vertex of $e_{i}$, then the Edge-Vertex Problem for $G$ is
equivalent to the ODD SET Problem for the red/blue bipartite graph
$G'=(R,B,E^*)$. But, this does not seem to help. In this section,
we will find all the graphs that have a solution to the
Edge-Vertex Problem, i.e., $G$ has even number of vertices. If $G$
has solutions, we give an algorithm to find a solution for $G$.
Some properties of the solutions and minimum solutions are
discussed. We get that a tree with even order has a unique
solution to the Edge-Vertex Problem, and therefore a unique
optimal solution to the Minimum Edge-Vertex Problem. Finally, by
polynomially transforming the Minimum Edge-Vertex Problem to the
Minimum Weight Perfect Matching Problem, we show that the Minimum
Edge-Vertex Problem for general graphs can be solved in polynomial time.\\

\noindent {\bf Definition 3.1.} For a graph $G$, if a spanning
subgraph $G'$ of $G$ has the property that the degree of every
vertex of $G'$ is odd in $G'$, then $G'$ is called an {\it odd
degree spanning subgraph}.\\

\noindent {\bf Theorem 3.1.} $G$ has a solution to the Edge-Vertex
Problem if and only if $G$ contains an odd degree spanning subgraph $G'$.\\

\noindent {\bf Proof. } If $G$ has a solution, then for all $v\in
V(G)$, the number of edges in the solution incident to the vertex
$v$ is odd. So, the subgraph $G'$ induced by the edges in the
solution is an odd degree spanning subgraph. The other round
is obvious. The proof is complete. \qed \\

The graphs that have a solution to the Edge-Vertex Problem are
given as follows.

\noindent {\bf Theorem 3.2.} A graph $G$ has a solution to the
Edge-Vertex Problem if and only if $G$ has even number of
vertices, i.e., $G$ is of even order.\\

\noindent {\bf Proof.} If $G$ has a solution to the Edge-Vertex
Problem, then from the above result, we know that $G$ contains a
odd degree spanning subgraph $G'$. Since the number of odd degree
vertices in any graph is even, we know that $G'$ has even number
of vertices. Since $G'$ is a spanning subgraph of $G$. so, $G$ has
even number of vertices.\\

Conversely, if $G$ has even number of vertices, we will show that
$G$ has a solution to the Edge-Vertex Problem. Since any solution
of a spanning tree of $G$ is also a solution of $G$, it is
sufficient to show that any tree $T$ of even order has a solution
to the Edge-Vertex Problem. By induction on $|T|=n$. If $|T|=2$,
the conclusion is obvious. Suppose that the conclusion is true for
all even order trees with the number of vertices less than $n$.
Then, for an even order tree $T$ with $|T|=n$, we can choose a
vertex $v$ from $T$ such that all but at most one neighbors of $v$
are not leaves of $T$. Such a vertex $v$ always exists by
observing that the vertices next to the end vertices of a longest
path of $T$ have the required property. We then distinguish
the following two cases:\\
[2mm] {\bf Case 1}. If all the neighbors of $v$ are leaves, then
$T=K_{1,n-1}$. It is easy to see that all the $n-1$ edges of
$K_{1,n-1}$ consist of a solution.\\
[2mm] {\bf Case 2}. Exactly one neighbor of $v$ is not a leaf of $T$.\\
[2mm] {\bf Subcase 2.1} Suppose the number of leaves adjacent to
$v$ is odd, say $l_{1},l_{2},\cdots,l_{2k+1}$. Let $T'$ be a
subtree of $T$ obtained by deleting $v$ and the leaves
$l_{1},l_{2},\cdots,l_{2k+1}$ from $T$. It is clear that $T'$ is
also an even order tree since $|T'|=n-(2k+1)-1=n-2(k+2)$. Then,
from the induction hypothesis, $T'$ contains a solution $E'$.
Therefore, $E'\cup
{\{vl_{1},vl_{2},\cdots,vl_{2k+1}\}}$ is a solution of $T$.\\
[2mm] {\bf Subcase 2.2} Suppose the number of leaves adjacent to
$v$ is even, say $l_{1},l_{2},\cdots,l_{2k}$. Let $T'$ be a
subtree of $T$ obtained by deleting the leaves
$l_{1},l_{2},\cdots,l_{2k}$ from $T$. It is clear that $T'$ is
also a tree of even order with $|T'|=n-2k$, and therefore contains
a solution $E'$ by induction hypothesis. We then obtain a solution
$E'\cup {\{vl_{1},vl_{2},\cdots,vl_{2k}\}}$ of $T$. The proof
is complete. \qed \\

Actually, the proof of the above theorem suggests an $O(n^2)$ or
$O(|E(G)|$ time algorithm for finding a solution to the
Edge-Vertex Problem for a given graph with even order, which is
stated as follows: \\

\noindent {\bf Algorithm to the Edge-Vertex Problem:} First, find
any spanning tree $T$ of $G$. Then, find a solution to the
Edge-Vertex Problem in $T$, by the method mentioned in the above
proof. All these can be done by the well-known BFS or DFS method.\\

\noindent {\bf Remark 3.1.} The complexity of the above algorithm
for a tree is obviously linear time, since each edge of a tree is
scanned at most once.\\

Before giving the algorithm for finding an optimal solution to the
Minimum Edge-Vertex Problem for a general graph, we show some
properties on optimal solutions.\\

\noindent {\bf Theorem 3.3.} If $G$ is a connected graph with even
number $n$ of vertices, then any optimal solution $S_{opt}$ to the
Minimum Edge-Vertex Problem for $G$ satisfies that $n/2\leq
|S_{opt}|\leq n-1$. Both the lower bound and the upper bound are
best possible.\\

\noindent {\bf Proof. } First, find a spanning tree $T$ of $G$,
and then we can get a solution to the Edge-Vertex Problem for $G$
from such a solution for $T$. So, $|S_{opt}|\leq n-1$. Since every
vertex of $G$ must be covered by $S_{opt}$ and every edge covers
at most two vertices, we get that $n/2 \leq |S_{opt}|$. To see
that the lower bound is best possible, we consider a graph $G$
with a perfect matching $M$. Then, $M$ is an optimal solution to
the Minimum Edge-Vertex Problem for $G$ and $|M|=n/2$. So, the
lower bound can be reached by many graphs. To see that the upper
bound is best possible, we consider the graph $G=K_{1,n-1}$. It is
not hard to see that the unique solution to the Edge-vertex
Problem for the graph $G$ is composed of all the edges of $G$. So,
the upper bound can be reached by the graph $G=K_{1,n-1}$. \qed \\

Actually, we can construct many other examples, as long as the
graph is a tree such that its every vertex has an odd degree.\\

The following result is immediate.\\

\noindent {\bf Theorem 3.4.} For any graph $G$ with even order
$n$, there is an $O(n^2)$ time algorithm for the Minimum
Edge-Vertex Problem of $G$ that produces a $2(1-1/n)$-approximation solution.\\

\noindent {\bf Proof. } First, find a spanning tree $T$ of $G$,
and then find a solution $S$ of $G$ from $T$. Since each of them
uses at most $O(n^2)$ time, the total complexity is $O(n^2)$.
Since the solution $S$ satisfies that\\
$$n/2\leq |S_{opt}| \leq|S| \leq n-1,$$
we have that
$$|S|/|S_{opt}|\leq(n-1)/(n/2)=2(1-1/n),$$\\
i.e., $S$ is a $2(1-1/n)$-approximation solution. \qed \\

\noindent {\bf Theorem 3.5.} The subgraph $G'$ induced by any
optimal solution to the Minimum Edge-Vertex Problem does not
contain any cycle. Moreover, each connected component of $G'$
induces a tree in $G$.\\

\noindent {\bf Proof. } Since every vertex of $G$ must be covered
by the solution, $G'$ must be a spanning subgraph of $G$. Next, we
want to show that there is no cycle in $G'$. If there is a cycle
$C$ in $G'$, then the edges of $C$ must belong to the optimal
solution. On the other hand, obviously $E(G')-E(C)$ is a solution
to the Edge-Vertex Problem for $G$, which contradicts to that
$E(G')$ is an optimal solution to the Minimum Edge-Vertex Problem
for $G$.\\

For the second statement of the theorem, suppose that a connected
component $H$ of $G'$ induces a subgraph that has a cycle. Since
$G'$ does not have any cycle and nor does $H$, there must be an
edge $e$ of $G$ such that there is a cycle $C$ in $H+e$. Then,
$(E(G')-E(C))\bigcup \{e\}$ is a solution to the Edge-Vertex
problem for $G$ with fewer edges than $E(G')$ since $G$ is simple
and so $|E(C)|\geq 3$, a contradiction to that $E(G')$ is optimal.
The proof is complete. \qed \\

From the above result, we get that the Minimum Edge-Vertex Problem
is equivalent to the following problem:\\

\noindent {\it The Problem of Odd Degree Induced Spanning Forest
with Minimum Number of Edges:}  In a given graph $G$, find an
induced spanning forest with Minimum number of edges such that
each of its vertex has an odd degree.\\

From the relation that if the number of components in a forest is
$r$ then the number of edges in the forest is $n-r$, we get that
the above problem is equivalent to the following problem:\\

\noindent {\it The Problem of Odd Degree Induced Spanning Forest
with Maximal Number of Components:} In a given graph $G$, find an
induced spanning forest with maximal number of components such
that each of its vertex has an odd degree. Or, in other words,
partition the set of vertices of $G$ into as many parts as
possible such that each part induces an odd degree tree in $G$.\\

Both the above two problems are interesting optimization problems.
But, it seems there is no help to solving our original Minimum
Edge-Vertex Problem. We will show that the Minimum Edge-Vertex
Problem can be solved in polynomial time. As a consequence, the
above two problems can also be solved in polynomial time.\\

The following facts are obvious: The order of every component in
the spanning forest of $G$ induced by an optimal solution is even,
and the number of components is at most $n/2$. The upper bound is
best possible, which can be reached by any graph with perfect
matchings. So, if a graph $G$ has a perfect matching, then a
solution $S$ is optimal if and only if $S$ is a perfect matching
of $G$.\\

There are many classes of graphs whose optimal solutions are
exactly their perfect matchings, for examples, $2$-edge-connected
and $3$-regular graphs have perfect matchings, $k$-regular
bipartite graphs with $k>0$ have perfect matchings, see [2]. Next
we find the following surprising result whose proof is very
simple. By {\it claw-free} in a graph $G$, we mean that $G$ does
not contain an induced $K_{1,3}$ as a subgraph, and we simply call
the graph a {\it claw-free
graph}. \\

\noindent {\bf Theorem 3.6.} A claw-free graph $G$ has a perfect
matching if and only if $G$ is of even order. \\

\noindent {\bf Proof. } If $G$ has a perfect matching, then
obviously $G$ is of even order. Conversely, if $G$ is of even
order, then $G$ has a solution to the Edge-Vertex Problem. Suppose
$S$ is an optimal solution to the Minimum Edge-Vertex Problem for
$G$. Then, each component of the subgraph of $G$ induced by $S$ is
an induced tree in $G$. Since $G$ is claw-free, we claim that each
such component is a $K_2$. Otherwise, some such component has a
path of length 2, say $uvw$. Since $v$ is of odd degree in the
induced subgraph, there must be a vertex $x$ different from $u$
and $w$ such that $vx$ is an edge in $S$. Then, from Theorem 3.5,
$\{u,v,w,x\}$ must induce a claw of $G$, a contradiction.
\qed \\

\noindent {\bf Remark 3.2.} Note that the line graph of any graph
is claw-free. So, a line graph with even order has perfect
matchings, and therefore the optimal solutions to the Minimum
Edge-Vertex Problem for a line graph are exactly the perfect
matchings of the line graph.\\

The relation between the optimal solutions and the perfect
matchings in a graph with perfect matchings suggests us the
following Berge-type results, the proofs are omitted.\\

\noindent {\bf Theorem 3.7.} For any two solutions $S$ and $S'$ to
the Edge-Vertex Problem for a graph $G$, the symmetric difference
$S\bigoplus S'$ induces an Eulerian subgraph of $G$.\\

\noindent {\bf Theorem 3.8.} Let $S$ be a solution to the
Edge-Vertex Problem for a graph $G$. Then, $S$ is optimal to the
corresponding minimum problem if and only if for any cycle $C$ of
$G$ we have that $|S\bigcap C|\leq |\overline{S}\bigcap C|$, where
$\overline{S}=E(G)-S$.\\

As a consequence, we get the following result.\\

\noindent {\bf Theorem 3.9.}  Suppose that $T$ is a tree with even
number of vertices. Then, $T$ has a unique solution to the
Edge-Vertex Problem, and therefore a unique optimal solution to
the Minimum Edge-Vertex Problem.\\

\noindent {\bf Proof. } If $H$ and $K$ are two different solutions
for $T$, then $H\bigoplus K$ is a nonempty Eulerian subgraph of
$T$. So, there is a cycle in the tree $T$, a contradiction. The
proof is complete. \qed \\

From Theorem 3.9, we can get a unique solution from a spanning
tree $T$ of a connected $G$ with even vertices. However, there is
no one to one correspondence between the solutions and the
spanning trees of a graph. Two different spanning trees may give
the same solution for $G$.\\

All the above discussions help us to understand the structures of
a (optimal) solution. Next we give the main result of this section.\\

\noindent {\bf Theorem 3.10.} The Minimum Edge-Vertex Problem can
be solved in polynomial time for general graphs $G$.\\

\noindent {\bf Proof.} It is easy to check if a graph $G$ has a
solution to the Edge-Vertex Problem can be done in $O(n)$, simply
by checking the order of $G$. To show that to find an optimal
solution to the Minimum Edge-Vertex Problem can be done in
polynomial time, we will polynomially transform the problem for a
given graph $G$ to the Minimum Weight Perfect Matching Problem for
another weighted graph $G^*$. Then, from [9] the later problem can
be solved in polynomial time, and so can be done for
the former problem. The proof is divided into the following steps.\\

\noindent {\bf Step 1.} Construct the graph $G^*$ from $G$.\\

Suppose that we are given a graph $G=(V(G),E(G))$ with even order
such that $V(G)=\{v_1,v_2,\cdots,v_n\}$. From our Algorithm to the
Edge-Vertex Problem, we can obtain a solution $S$ to the
Edge-Vertex Problem for $G$ in $O(n^2)$ time. Consider the
subgraph $G_0$ of $G$ induced by the $S$. We know that $G_0$ is an
odd degree spanning subgraph of $G$. For every vertex $v$ of $G$,
assume that the neighborhood of $v_i$ in $G$ is
$N_i=\{u_1,u_2,\cdots,u_d\}$, and the neighborhood of $v_i$ in
$G_0$ is $\{u_1,u_2,\cdots,u_{2k-1}\}$, where $k\geq 1$ and
$2k-1\leq d$.
We distinguish the following cases.\\

\noindent {\bf Case 1.1.} $d$ is odd, say $d=2r-1$.\\

Obviously, $r\geq k$. Let $V_i^1=\{v_i^{11},v_i^{12},\cdots,
v_i^{1d}\}$ and $V_i^2=\{v_i^{21},v_i^{22},\cdots,v_i^{2d}\}$ be
two new sets of vertices corresponding to the vertex $v_i$.
Completely connect the vertices in $V_i^1$ to form a complete
graph with $d$ vertices, and assign each of the edges with weight
0. Then, connect $v_i^{1j}$ and $v_i^{2j}$ with an edge for every
$j=1,2,\cdots,d$, and assign each of these edges with weight 1.
Finally, connect $v_i^{2j}$ and $u_j$ with an edge for every
$j=1,2,\cdots,d$, and assign each of these edges with weight 0.
The edge $v_i^{2j}u_j$ can be imagined to be equal to the edge
$v_iu_j$ for every $j=1,2,\cdots,d$.\\

\noindent {\bf Case 1.2.} $d$ is even, say $d=2r$.\\

Obviously, $r\geq k$. Let $V_i^1$ and $V_i^2$ be the same sets as
in Case 1.1. The edges among the vertices in the three sets
$V_i^1$, $V_i^2$ and $N_i$ are also the same as in Case 1.1.
However, we add one more new vertex $w_i$ for $v_i$. Then, connect
$w_i$ and $v_i^{1j}$ with an edge for every $j=1,2,\cdots,d$, and
assign each of these edge with weight 0. So, $V_i^1\bigcup
\{w_i\}$ forms a complete graph with $d+1$ vertices.\\

Let $W=\{w_i | d_G(v_i)=0 \mod 2\}$. Then, the new graph $G^*$ has
$W\bigcup (\bigcup_{i=1}^{n}V_i^1)\bigcup (\bigcup_{i=1}^n V_i^2)$
as its vertex set. The edge set of $G^*$ consists of all the edges
in Cases 1.1 and 1.2 with a small change for the kind of edges
$v_i^{2j}u_j$, because then $u_j$ is also blown up into two sets
of vertices (plus one more vertex if $d_G(u_j)$ is even). In
general, we call the vertices $v_i^{2j}$ for $j=1,2,\cdots,d$ the
representatives of the vertex $v_i$ in $G^*$. If $v_iv_j$ is an
edge of $G$, then $v_j$ is a neighbor of $v_i$ and $v_i$ is a
neighbor of $v_j$. Suppose that $v_j$ is the $p$-th neighbor of
$v_i$ and $v_i$ is the $q$-th neighbor of $v_j$. Then we use the
edge $v_i^{2p}v_j^{2q}$ in $G^*$ to substitute the edge
$v_i^{2p}v_j$ as in Cases 1.1 and 1.2, which represents the edge
$v_iv_j$ of $G$. Now the construction of the new graph $G^*$ is
ready. Obviously, the construction can be done in polynomial time
in $n$.\\

\noindent {\bf Step 2.} Show that $G^*$ has a perfect matching.\\

We construct a perfect matching $M$ with the following three
parts:\\

\noindent {\bf Part 2.1.} Recall that we have an odd spanning
subgraph $G_0$. For every edge $v_iv_j$ in $E(G)-E(G_0)$, we order
that the corresponding edge $v_i^{2p}v_j^{2q}$ in $G^*$ belongs to $M$.\\

Note that these matched even number of vertices in $V_i^2$ if
$v_i$ has an odd degree in $G$, and odd number of vertices in
$V_i^2$ if $v_i$ has an even degree in $G$, for every $i=1,2,\cdots,n$.\\

\noindent {\bf Part 2.2.} For each odd degree vertex $v_i$, we
match the even number of vertices $v_i^{1j}$ in $V_i^1$
corresponding to the matched vertices $v_i^{2j}$ in $V_i^2$ by any
perfect matching among them. Then, match the remaining unmatched
vertices in $V_i^1$ and $V_i^2$ by $v_i^{1j}v_i^{2j}$.\\

Note that all the vertices blown up at $v_i$ are matched now.\\

\noindent {\bf Part 2.3.} For each even degree vertex $v_i$, we
match the odd number of vertices $v_i^{1j}$ in $V_i^1$
corresponding to the matched vertices $v_i^{2j}$ in $V_i^2$ by a
perfect matching among them plus the extra one vertex $w_i$. The
rest can be done in the same way as in Part 2.2.\\

Obviously, the above constructed is a perfect matching for the
graph $G^*$.\\

\noindent {\bf Step 3.} Show that every solution to the
Edge-Vertex Problem for the original graph $G$, i.e., the set of
edges of every odd degree spanning subgraph of $G$ corresponds to
a perfect matching of $G^*$.\\

The corresponding perfect matching $G^*$ can be constructed in the
same way as in Step 2.\\

\noindent {\bf Step 4.} Denote the set of edges $v_i^{2p}v_j^{2q}$
in $G^*$ corresponding to the edges $v_iv_j$ in $G$ by $E^*(G)$.
Then there is a one to one correspondence $\varphi$ from $E^*(G)$
onto $E(G)$. Show that every perfect matching $M$ of $G^*$
corresponds, under $\varphi |_{\overline{M}\bigcap E^*(G)}$, to
the set of edges of an odd degree spanning subgraph of $G$, i.e.,
a solution to the Edge-Vertex Problem for the graph $G$.\\

We only need to show that $M$ always matches odd number of edges
between $V_i^1$ and $V_i^2$ for any vertex $v_i$, since then for
each matched edge between $V_i^1$ and $V_i^2$ there is exactly one
unmatched edge between $V_i^2$ and $N_i$, which corresponds to the
edges in an odd degree spanning subgraph of $G$. We distinguish
two cases.\\

\noindent {\bf Case 4.1.} $v_i$ is an odd degree vertex of $G$.\\

Since $M$ is a perfect matching, all the vertices in $V_i^1$ must
be matched. In any case, there are even number of vertices matched
among $V_i^1$. Then, the rest odd number of vertices in $V_i^1$
unmatched among $V_i^1$ must matched with odd number of vertices
in $V_i^2$. Our required is proved.\\

\noindent {\bf Case 4.2.} $v_i$ is an even degree vertex in $G$.\\

Since $M$ is a perfect matching, all the vertices in $V_i^1\bigcup
\{w_i\}$ must be matched. Note that the vertex $w_i$ must be
matched with a vertex in $V_i^1$. In any case, the rest odd number
of vertices in $V_i^1$ can match even number of them among
themselves. So, the number of vertices of $V_i^1$ matched with the
vertices in $V_i^2$ is odd. Again, the required is proved.\\

\noindent {\bf Step 5.} Note that $G^*$ is a weighted graph with
perfect matchings. Each of the edges between $V_i^1$ and $V_i^2$
is assigned weight 1, while each of the other edges is assigned
weight 0. Show that the weight of a minimum weight perfect
matching $M$ of $G^*$ is equal to the number of edges of an
optimal solution to the Minimum Edge-Vertex Problem for $G$, or
the number of edges of an odd degree spanning subgraph of $G$ with
minimum number of edges.\\

From the proof of Step 4, a minimum weight perfect matching $M$ of
$G^*$ corresponds to an odd degree spanning subgraph of $G$ with
the number of edges equal to the weight of $M$. So, the weight of
a minimum weight perfect matching is at least the number of edges
of an odd degree spanning subgraph of $G$ with Minimum number of
edges. Conversely, from Step 2, an odd degree spanning subgraph
$G_0$ of $G$ with Minimum number of edges can be used to construct
a perfect matching $M$ of $G^*$. It is easy to see that the number
of edges in $G_0$ is equal to the weight of $M$. So, the number of
edges in an odd degree spanning subgraph of $G$ with Minimum
number of edges is at least the weight of a minimum weight perfect
matching of $G^*$. Therefore, the minimum weight and the Minimum
number must be equal.\\

Up to now, we have proved that the Minimum Edge-Vertex Problem for
a given graph $G$ is equivalent to the Minimum Weight Perfect
Matching Problem for a new graph $G^*$ polynomially constructed
from $G$. From [9] we know that the Minimum Weight Perfect
Matching Problem can be solved in polynomial time. So, the Minimum
Edge-Vertex Problem can also be solved in polynomial time. The
proof is now complete. \qed \\

\noindent {\bf Remark 3.3.} The above theorem tells us that three
other problems can also be solved in polynomial time, i.e., the
Problem of Odd Degree Induced Spanning Forest with Minimum Number
of Edges, the Problem of Odd Degree Induced Spanning Forest with
Maximal Number of Components, and the Odd Set Problem [5] for a
red/blue bipartite graph $(R,B,E)$ such that each vertex in $R$
has degree 2.\\

\section{Concluding Remarks}

As we mentioned before, it is easy to see that the Edge-Edge
Problem for a given graph $G$ is equivalent to the Vertex-Vertex
Problem for the line graph $L(G)$ of $G$. So, the existence of
solutions to the Edge-Edge Problem for general graphs is solved,
and to find solutions can also be done in polynomial time.
However, for the Minimum Edge-Edge Problem, it seems no help for
finding an optimal solution efficiently for this special kind of
graphs, i.e., line graphs, which have many particular properties,
for example, claw-free. We shall leave the minimum (optimization)
problem for further study.\\
[6mm] \noindent {\bf Acknowledgement.} The authors would like to
thank Prof. G. Woeginger of Twente University, The Netherlands,
for helpful discussions.\\
[5mm]

\end{document}